\documentclass[11pt]{article}   
 \pdfoutput=1 
 \usepackage{graphicx}
\usepackage{amsthm,amsmath,amssymb}   
\usepackage{color}
\usepackage{url}
\usepackage{enumerate}
\usepackage{epsfig,subfigure}
\newtheorem{theorem}{Theorem}[section]

\newtheorem{definition}[theorem]{Definition}

\theoremstyle{remark}
\newtheorem{remark}[theorem]{Remark}

\title{\LARGE Dynamic constructions of hyperbolisms\\ of plane curves:\\
an automated exploration of geometric loci}

\date{\today}

\begin{document}

\maketitle
\begin{center}
{\bf  Thierry Dana-Picard}\\[1mm]
Jerusalem College of Technology and Jerusalem Michlala College \\ 
Jerusalem, Israel \\[1mm]
{\tt  ndp@jct.ac.il}
\end{center}

\setcounter{page}{1}

\vspace*{1mm}

\begin{abstract}
{\it Hyperbolism of a given curve with respect to a point and a line is an interesting construct, a special kind of geometric locus, not frequent in the literature.  While networking between two different kinds of mathematical software, we  explore various cases, involving quartics, among them the so-called K\"ulp quartic and topologically equivalent curves, and also an example with a sextic and a curve of degree 12. By a similar but different way, we derive a new construction of a lemniscate of Gerono. First, parametric equations are derived for the curve, then we perform implicitization Gr\"obner bases packages and using elimination. The polynomial equation which is obtained enables to check irreducibility of the constructed curve.}
\end{abstract}
\vspace*{1mm}

\section{Introduction}
\label{intro}
 
 \subsection{The needed dialog between two kinds of mathematical software}
 Prior to the development of Computer Algebra Systems (CAS) and Dynamic Geometry Software (DGS), mechanical devices were used to draw specific curves \cite{historical mechanisms}; spirographs were among the most popular and enabled to draw epitrochoids and hypotrochoids\footnote{See \url{https://mathcurve.com/courbes2d.gb/epitrochoid/epitrochoid.shtml} and \url{https://mathcurve.com/courbes2d.gb/hypotrochoid/hypotrochoid.shtml}.} and other related curves. Here, instead of mechanical devices, we use software, both a DGS (GeoGebra) and a CAS (Maple). The first one enables to construct and explore the curves, their shape and topology, sometimes providing polynomial equations but not always. We use then the CAS in order to perform algebraic computations and derive polynomial equations for the curves under study.  Afterwards it is possible to copy-paste the formulas to the DGS and to analyze the curves with its dynamical features. This CAS-DGS collaboration provides a useful environment for the exploration of new constructs, as in \cite{DP-K-Y}. Such  a dialog between the two kinds of software has been used for example in \cite{DPK-networking,maltese},and for years has been wished to be more automatic \cite{eugenio}. 
 
In this work,  the curve are determined first by parametric equations. Implicitization is important, as a polynomial presentation may provide information whether the curve is irreducible or not (see \cite{fischer}, chap. 1). The algebraic part of the work  is based on the theory of Gr\"obner bases and on Elimination; see \cite{cox,sendra winkler}. We may refer also to \cite{adams and loustaunau}, in particular for the non familiar reader who can see how to work out elementary examples by running the algorithms "by hand". Note that Maple has an \textbf{implicitize} routine in its package \emph{algcurves}, based on  \cite{corless kotsireas et al}. Irreducibility is also checked with automated methods; see Section \ref{section circle and ellipse}.

Mathematical objects cannot be grasped with hands, and are approached using numerous registers of representations \cite{duval}; the classical registers for plane curves are graphical, numerical and algebraic. Parametric representations and implicit representations as two subregisters of the algebraic one. The study is made rich and efficient by switching between registers, but switching from parametric to implicit and from implicit to parametric are non trivial tasks; see \cite{implicitizaton gao et al,sederberg and goldman,sendra winkler,cox}. In some cases, the switch is impossible.  

Switching between different representations of plane curves is an important issue, with numerous applications in computer aided design and other fields. In \cite{wang}, Wang emphasizes the role in computer aided geometric design and modeling. He develops "an extremely simple method that converts the rational parametric equations for any curve or surface into an implicit equation" (we recommend also the vast bibliography there in the paper) . He uses Gr\"obner bases, resultants, etc. In our work here, we transform the obtained parametric equations into parametric rational presentations, then into polynomial equations. We use Maple's \emph{PolynomialIdeals} package and elimination to derive implicit equations. It is often easier to analyze the topology of a curve using an implicit presentation than a parametric presentation. Anyway, both enable the study and classification of singular points.

 \subsection{Plane curves defined as geometric loci}
\label{subsection geom loci and hyperbolisms}
 Plane algebraic curves are a classical topic, to which numerous books have been devoted, such as \cite{yates}. Websites are devoted to curves (and surfaces) such as Mathcurve (\url{http://mathcurve.com}). Full catalogues of curves of degree 2,3 and 4 exist, and partial catalogues for degree 6.  Some of them are constructed as  geometric loci. The bifocal definition of ellipses and hyperbolas is generally the first example met by students. Cassini ovals (also called spiric curves) and Cayley ovals are more advanced examples. Recently, some octic curves (curves defined by polynomials of degree 8)  have been described as geometric loci in relation with a classical theorem of plane geometry, namely Thales second theorem \cite{DPR, DP-R}.  The present paper shows a bunch of curves, more or less classical, appearing as geometric loci by  construction of \emph{hyperbolisms}.
 
Our concern is the construction and study of plane curves as hyperbolism of classical curves, according to the definition in the Mathcurve website \cite{hyperbolism mathcurve}:
\begin{definition}
\label{def hyperbolism}
The hyperbolism of a curve  $\Gamma_0$ with respect to a point $O$ and a line $(D)$ is the curve $\Gamma$, locus of the point $M$ defined as follows: given a point $M_0$ on $\Gamma_0$, the line $(OM_0)$ cuts $(D)$ at $P$; $M$ is the projection of $P$ on the line parallel to $(D)$ passing by $M_0$.
\end{definition}
 \begin{remark}
 Analytically, if the line $(D)$ is given  by the equation $ x = a$, the transformation of $\Gamma_0$  into $\Gamma$ can be written $(x,y) \rightarrow \left( x, a \frac {y}{x} \right)$; it is quadratic, so an algebraic curve of degree $n$ is transformed into an algebraic curve of degree at most $2n$. The 2nd coordinate $a \frac {y}{x}$ makes the connection with hyperbolas, whence their name \emph{hyperbolism}.
 \end{remark}

 By definition, hyperbolisms of curves are a subtopic of geometric loci, which are a classical topic, from middle school to university. The last decades have seen numerous works devoted to automated methods in Geometry, such as  \cite{taxonomy,botana-abanades,blazek-pech} and  \cite{ART}, where loci are one of the main topics for which automated methods have been developed, and implemented in GeoGebra-Discovery\footnote{A freely downloadable companion to GeoGebra, available from\url{https://github.com/kovzol/geogebra-discovery}}.  A large number of versions of automated commands exist, we use here only a few of them, mostly  GeoGebra's \textbf{Locus($<$Point Creating Locus$>$,$<$Point$>$)}  and \textbf{Locus($<$Point Creating Locus$>$,$<$Point$>$)}. The place holders are called respectively the \emph{Tracer} and the \emph{Mover} in the above mentioned papers. The automated command provides a plot of the curve, sometimes also an implicit equation. For this it uses numerical methods. In the companion package GeoGebra-Discovery, this has been supplemented by symbolic algorithms yielding more precise answers. For the algebraic work, we switched to  the Maple software, especially for implicitization, but not only.   
  
 The notion of a hyperbolism has been presented to a small group of in-service mathematics teachers, learning towards an advanced degree M.Ed. These teachers had previous knowledge including the perpendicular bisector of a segment, a circle, conics, etc. as geometric loci, but almost neither CAS nor DGS  literacy. The outcome of their work was double: the development of new perspectives in geometry and acquisition of technological skills. They had neither true and the general atmosphere of the course was to  show new mathematical topics in a technology-rich environment and develop technological skills. The feedback was very positive, and true curiosity at work.
  
 \section{First easy examples}
 \label{section circle and ellipse}
 \subsection{Hyperbolism of a circle centered at the point $O$ and the line is tangent to the circle}
 \label{subsection circle 1}
\begin{figure}[htb]
\begin{center}
\includegraphics[width=11cm]{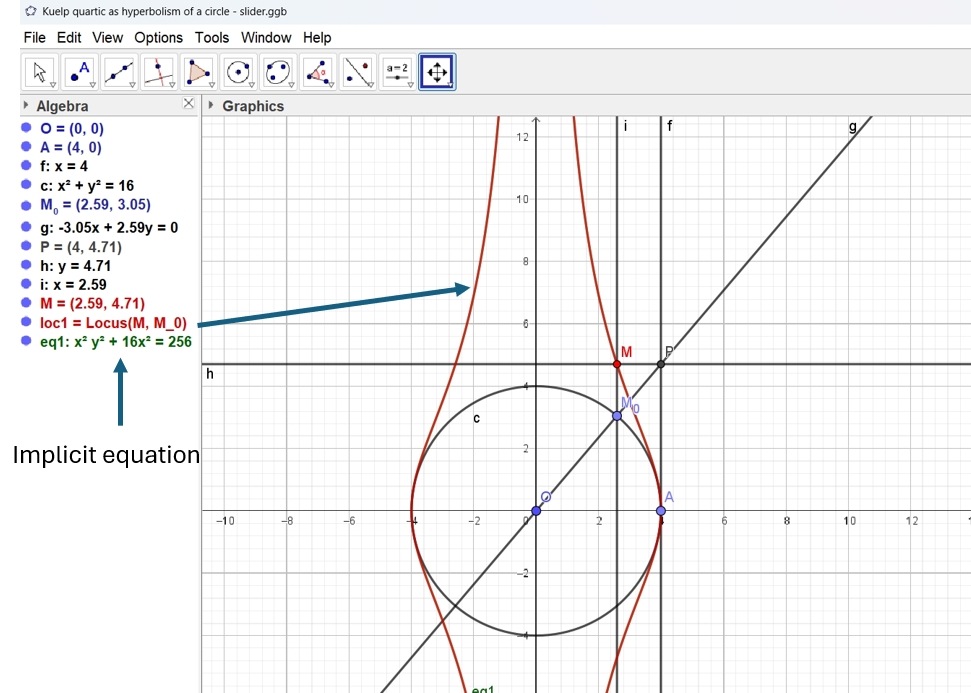}
\caption{Hyperbolism of a circle with respect to its center and a tangent}
\label{fig hyperbolism of a circle 1}
\end{center}
\end{figure}

 We consider the circle $\mathcal{C}$ centered at the origin $O$ with radius $r$. The line $D$ has equation $x=r$ and is thus tangent to the circle $\mathcal{C}$. A point $M_0 \in \mathcal{C}$ is given by $(x,y)=(r \cos t, r \sin t)$ where  $t \in [0,2\pi ]$. The line $OM_0$ has thus equation $y=x \tan t$ and intersects the line $D$ at $P(r \cos t,r \tan t)$. It follows that the point $M$ has coordinates 
 \begin{equation}
 \label{param eq hyperbolism of a circle 1}
 (x,y)=(r \cos t, r \tan t)
 \end{equation}
 Equation (\ref{param eq hyperbolism of a circle 1}) is a parametric presentation of the geometric locus that we are looking for. Figure \ref{fig hyperbolism of a circle 1} shows a screenshot of a GeoGebra session for this question. The requested curve has been obtained using the \textbf{LocusEquation($M,M_0$)} command. It could have been obtained also with the \textbf{LocusEquation($<$Point Creating Locus Line $>$,$<$ Slider $>$)}, after entering directly the parametrization (\ref{param eq hyperbolism of a circle 1}), but in this case the new construct is independent of what has been done previously.
The plots corresponding to the 2 last rows in the algebraic window overlap each other, therefore only one plot is viewed. Nevertheless, they are considered by the software as 2 different objects.

\begin{figure}[htb]
\begin{center}
 \includegraphics[width=11cm]{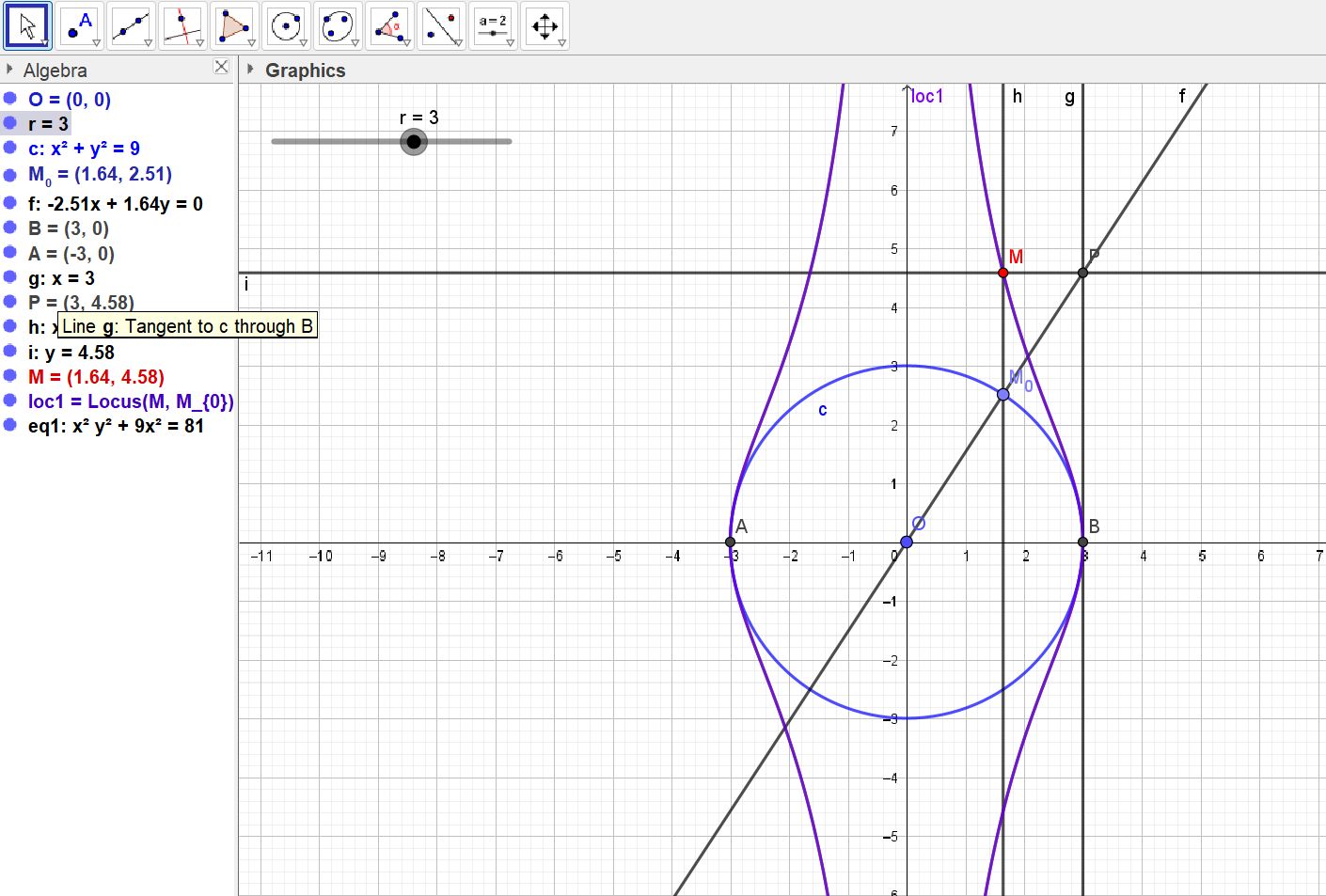}
\caption{Hyperbolism of a circle with respect to its center and a tangent - variable radius}
\label{fig hyperbolism of a circle 2}
\end{center}
\end{figure}

In this session, the radius may be changed, but the equation of the geometric locus remains of the same form. In Figure \ref{fig hyperbolism of a circle 1}, $r=4$ and the geometric locus has equation
\begin{equation}
\label{eq kulp quartic 1}
x^2y^2+16x^2=256,
\end{equation}
which is identified as the equation of a so-called \emph{K\"ulp quartic}\footnote{A curve studied in 1878 by K\"ulp; see \url{https://mathcurve.com/courbes2d/kulp/kulp.shtml}}.
A GeoGebra session\footnote{\url{https://www.geogebra.org/m/md4f8aaa}} may help to visualize the construct for various values of the radius. Slightly different ways to construct the hyperbolism may be chosen, and for some of them an implicit equation is not available. The user has to modify his protocol to have it "more geometric" (v.i. Section \ref{section piriform}). Sometimes, GeoGebra-Discovery displays a message telling that the construction involves steps which are not supported by the command \textbf{LocusEquation}. In our applet, the tangent has not been constructed directly with its equation, as in Figure \ref{fig hyperbolism of a circle 1}, but using the \textbf{Tangent} command of the software, which is important as part of a geometric construct; see Figure \ref{fig hyperbolism of a circle 2}.

In order to derive from Equation (\ref{param eq hyperbolism of a circle 1}) a polynomial equation, we use the following substitution, as in \cite{DPK-networking}:
\begin{equation}
\label{rational substitution for trigo}
\forall t \in \mathbb{R}, \; \exists u \in \mathbb{R} \; \text{such that} \; 
\begin{cases}
\cos t = \frac{1-u^2}{1+u^2}\\
\sin t = \frac{2u}{1+u^2}
\end{cases}
\end{equation}
We apply the following Maple code:
\small
\begin{verbatim}
xk := r*cos(t); yk := r*sin(t)/cos(t);
xkrat := subs(cos(t) = (-u^2 + 1)/(u^2 + 1), xk);
ykrat := subs(cos(t) = (-u^2 + 1)/(u^2 + 1), 
         subs(sin(t) = 2*u/(u^2 + 1), yk));
p1 := x*denom(xkrat) - numer(xkrat);
p2 := y*denom(ykrat) - numer(ykrat);
J := <p1, p2>;
JE := EliminationIdeal(J, {r, x, y});
\end{verbatim}
The output is an ideal with a unique generator, providing the following equation (for any real $r$):
  \begin{equation}
 \label{eq Kulp quartic}
r^2x^2 + x^2y^2-r^4=0.
 \end{equation}
Note the dependence of the coefficients in Equation (\ref{eq Kulp quartic}) on the radius of the circle.
 Using Maple's command \textbf{evala(AFactor(...)}, we check that the left hand side in  Equation (\ref{eq Kulp quartic}) is an irreducible polynomial. This means that the 2 components of the curve described by this equation cannot be distinguished by algebraic means, i.e. they are not 2 distinct components of the curve (see \cite{fischer}, p.17-18).

 \subsection{Hyperbolism of an ellipse}
 \label{subsection Hyperbolism of an ellipse}
 We generalize slightly the situation of subsection \ref{subsection circle 1} and consider now an ellipse $\mathcal{C}$ whose equation is $\frac{x^2}{a^2}+\frac{y^2}{b^2}=1$, where $a$ and $b$ are positive parameters\footnote{GeoGebra applets are available at  \url{https://www.geogebra.org/m/eqywkwcw} and \url{https://www.geogebra.org/m/tfnrwaqh}}. Here too, the point $O$ is the origin. The line $D$ has equation $x=a$, i.e. is tangent to the ellipse $\mathcal{C}$ and parallel to the $y-$axis.
 A point $M_0 \in \mathcal{C}$ is given by $(x,y)=(a \cos t, b \sin t)$ where  $t \in [0,2\pi ]$. The line $OM_0$ has thus equation $y=x \frac{b}{a}\tan t$ and intersects the line $D$ at $P(a ,b \tan t)$. It follows that the point $M$ has coordinates 
 \begin{equation}
 \label{param eq hyperbolism of an ellipse}
 (x,y)=(a \cos t, b \tan t)
 \end{equation}
 \begin{figure}[htb]
\begin{center}
 \includegraphics[width=11cm]{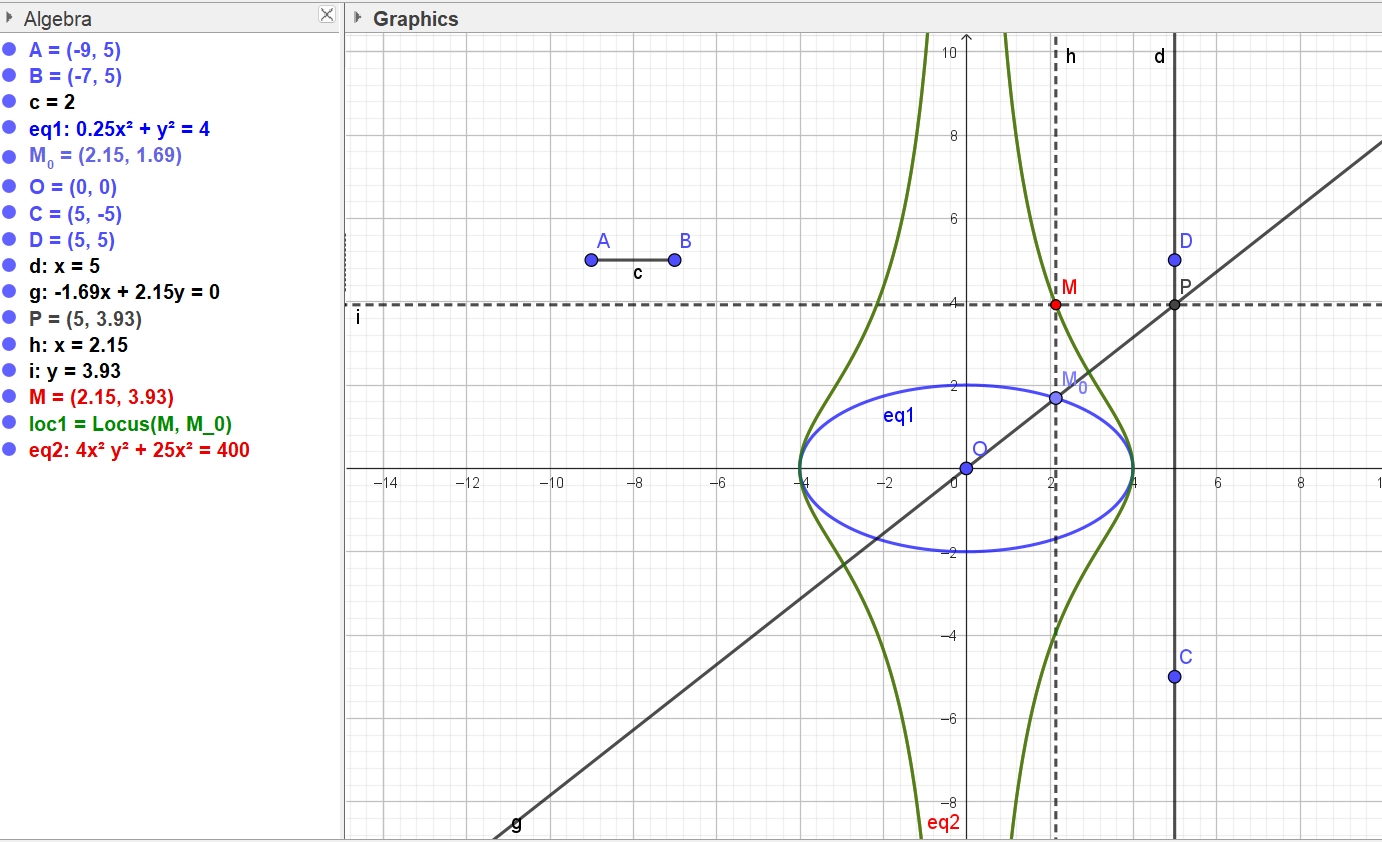}
 \caption{Hyperbolism of an ellipse}
\label{fig hyperbolism of an ellipse}
\end{center}
\end{figure}

\begin{remark}
\label{remark degree}
The segment $c$ in the upper left corner of Figure \ref{fig hyperbolism of an ellipse} comes instead of a slider for a  parameter which is involved in the implicit equation of the ellipse. This, in order to have the ellipse dependent on a geometric construct. Actually, the ellipse itself could have been constructed as a geometric loci, but in such a case GeoGebra's command \textbf{LocusEquation} may provide a plot and an implicit equation, but not the possibility to use the command \textbf{Point on Object}.
\end{remark}

\begin{figure}[htb]
\begin{center}
 \includegraphics[width=11cm]{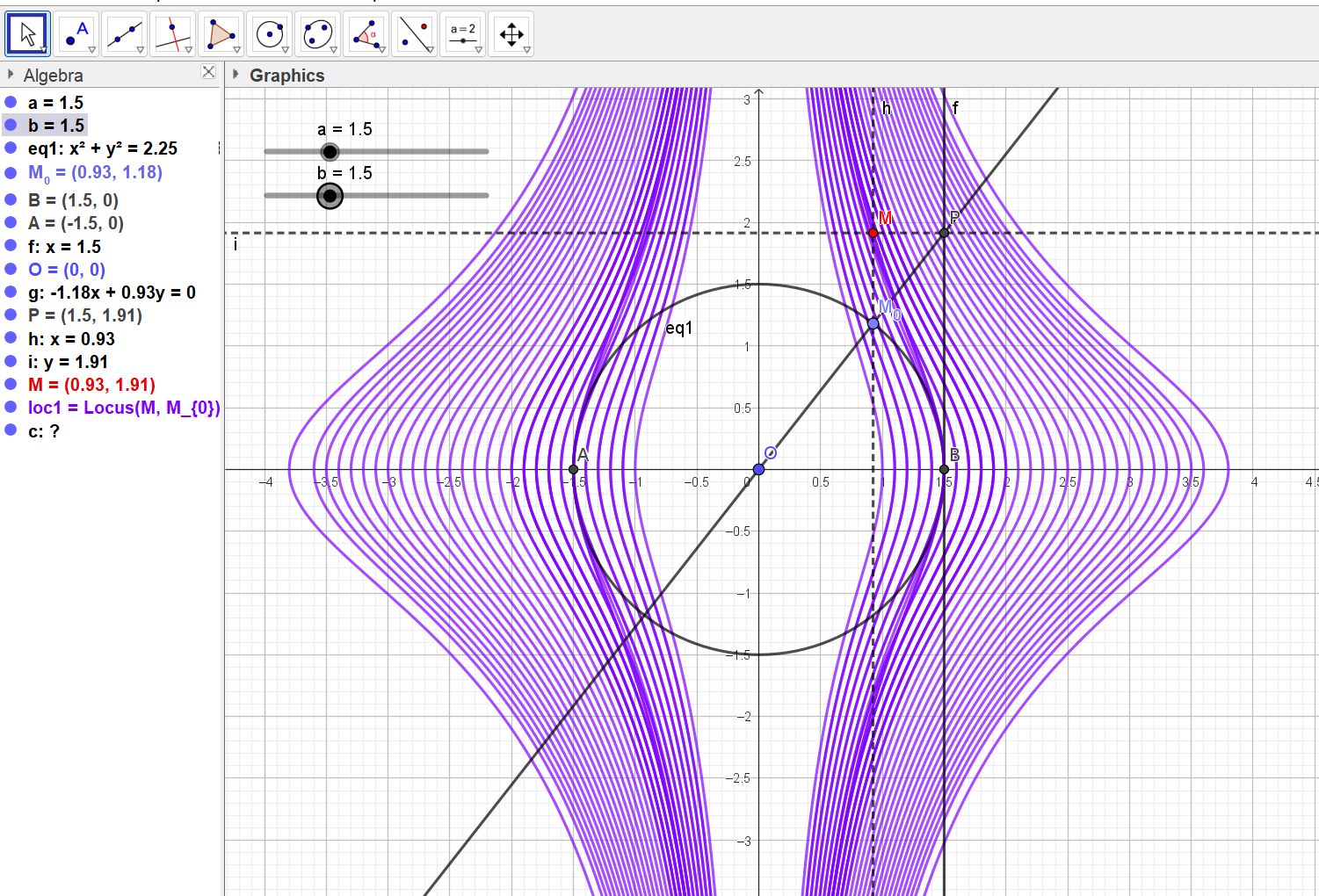}
 \caption{Showing that all the curves have the same topology}
\label{fig explore topology}
\end{center}
\end{figure}

We apply a Maple code similar to the code in the previous subsection\footnote{Maple's \textbf{implicitize} command did not always provide an answer.}. The output provides an equation for the hyperbolism of the ellipse:
\begin{equation}
\label{eq hyperbolism of an ellipse}
b^2x^2 + x^2y^2 -a^2b^2=0.
\end{equation}
Using the same tool as in previous section, we show that obtained curve is irreducible. Note that it is not a K\"ulp quartic (the coefficient of $x^2$ is not the square of the free coefficient). The equations are quite similar, but not identical. Figure \ref{fig explore topology} shows how a GeoGebra applet may help to understand that, for various values of the axes of the ellipse, the obtained curves have the same topology. The exploration leads to proving that two curves in this family are obtained from each other by an affinity whose axis is the $y-$axis and the direction is perpendicular to it. Moreover, this may be an opportunity to explore in class the similarities and the differences of these quartics, K\"ulp quartic, and also the Witch of Agnesi (which is a cubic).

\begin{remark}
Denote $F(x,y)= b^2x^2 + x^2y^2 -a^2b^2$. We have: $F_x(x,y)=2x(b^2+y)$ and $F_y(x,y)=2x^2y$. It is easy to show that the system of equations $F_x=F_y=0$ has no solution (actually $(0,0)$ solves the system, but the origin does not belongs to the curve). Therefore the curve has no singular point.
\end{remark}

\begin{remark}
Assuming that at inflexion points the curvature is equal to 0, it is possible to look for candidates using the following code:
\small
\begin{verbatim}
with(Student[VectorCalculus]);
Curvature(<a*cos(t), b*tan(t), t>, t);
simplify(%);
infl := solve(% = 0, t);
allvalues(infl[1]);
\end{verbatim}
Different issues may appear. First the meanings of $<...,...>$ in the two packages \emph{PolynomialIdeals} and \emph{VectorCalculus} are different, therefore we suggest to introduce the 2nd one only after the first one has been used. Second, the answer for general parameters $a$ and $b$ may be heavy; it may be wiser to use them with specific values. Finally, this provides values of the parameter $t$ which can correspond to points of inflexion, but not only. More verifications are needed.
\end{remark}

\begin{remark}
The given ellipse is the image of the circle whose equation is $x^2+y^2=a^2$ by the affine transformation $(x,y) \longmapsto \left( x, \frac {b}{a}y \right)$. It is easily proven that the quartic obtained here is the image  by the same transformation of the K\"ulp quartic found in the previous subsection.
\end{remark}

\section{Hyperbolism of a circle with respect to a line secant to the circle}
\subsection{The circle is centered at the origin}
Other quartics can be obtained as hyperbolisms of circles. 

\begin{figure}[htb]
\begin{center}
 \includegraphics[width=11cm]{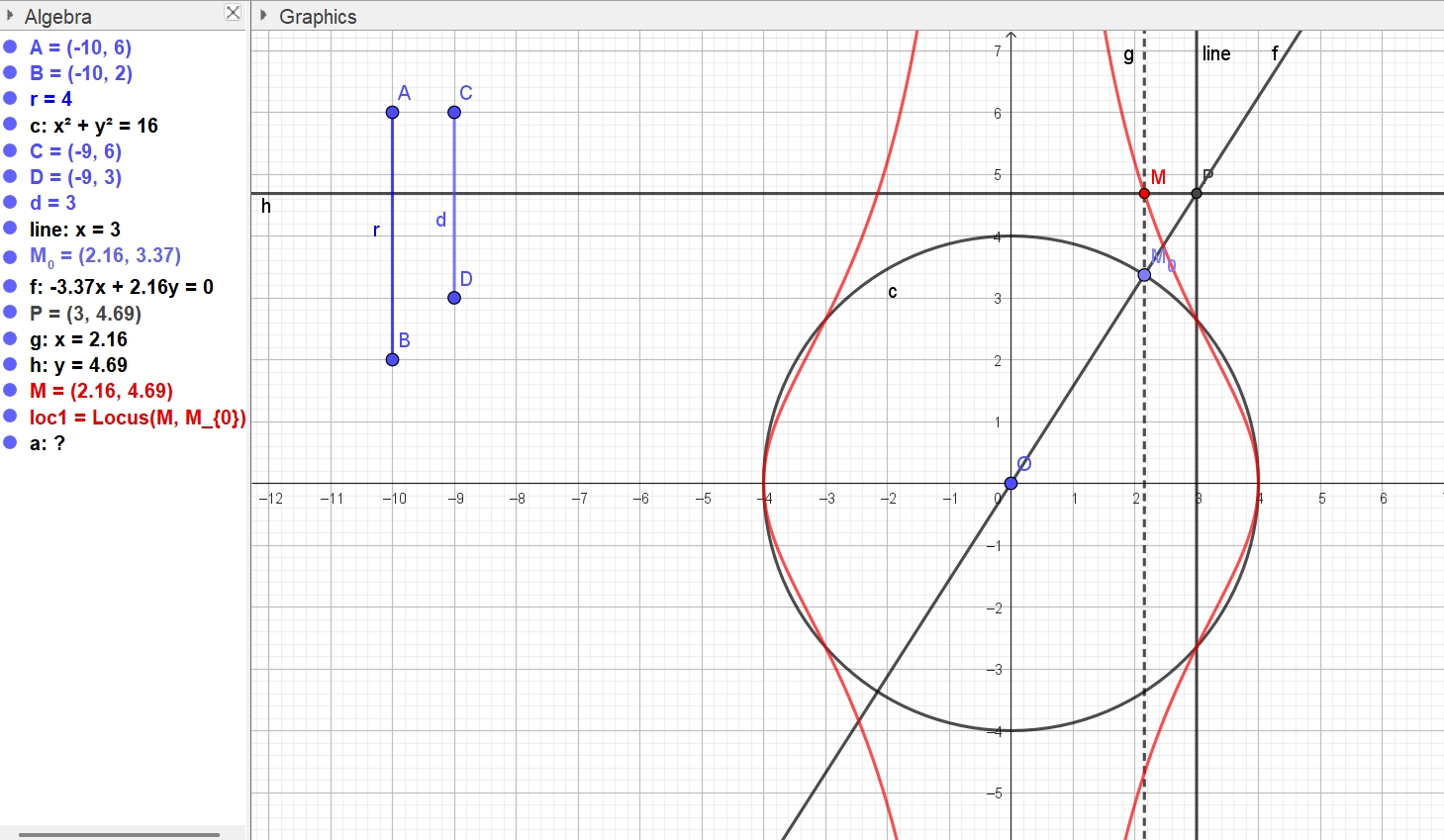}
 \caption{The line is not tangent to the circle}
\label{fig line not tangent to the circle}
\end{center}
\end{figure}

\begin{figure}[h]
\begin{center}
 \includegraphics[height=5cm]{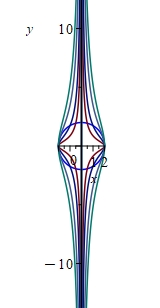}
 \caption{The line is not tangent to the circle - several plots}
\label{fig circle centered at origin - non tangent line- multiple plots}
\end{center}
\end{figure}

We consider a circle $\mathcal{C}$ centered at the origin with radius $r$ and the line $D$ with equation $x=b$, where $b \neq r$. Figure \ref{fig line not tangent to the circle}  is a screenshot of a GeoGebra applet\footnote{\url{https://www.geogebra.org/m/pwrx2uzf}; a related applet is available at \url{https://www.geogebra.org/m/x22etuhc}.}); the line $D$ is a secant to the circle $\mathcal{C}$.
We perform the same construction as in subsection \ref{subsection circle 1}, and obtain a plot of the geometric locus, but cannot obtain an implicit equation with the command \textbf{LocusEquation}. Therefore, we need to make the algebraic computation using the CAS. Note that, even when changing the value of $fr$ by changing the length of the corresponding segment AB, the implicit equation is not obtained.
\small
\begin{verbatim}
f := -r^2 + x^2 + y^2;
li := y = r*sin(t) + sin(t)*(x - r*cos(t))/cos(t);
yP := subs(x = d, rhs(li));
xM := r*cos(t);
yM := yP;
p1 := x - xM;
p2 := y - yM;
p1 := subs(cos(t) = (-u^2 + 1)/(u^2 + 1), p1);
p2 := subs(cos(t) = (-u^2 + 1)/(u^2 + 1), subs(sin(t) = 2*u/(u^2 + 1), p2));
p1 := simplify(p1*denom(p1));
p2 := simplify(p2*denom(p2));
J := <p1, p2>;
JE := EliminationIdeal(J, {d, r, x, y});
\end{verbatim}
\normalsize
The output provides an implicit equation for the hyperbolism. After simplification, it reads as follows:
\begin{equation}
\label{implicit eq circle centered at the origin - line parallel to y-axis}
 x^2y^2+ d^2x^2 -d^2r^2=0
\end{equation}
As expected, for $d=r$, we have Equation (\ref{eq Kulp quartic}).
Figure \ref{fig circle centered at origin - non tangent line- multiple plots} shows hyperbolisms of the same circle for  4 different lines whose respective equations are $x=1/2$, $x=1$, $x=2$ and $x=3$, from the innermost to the outermost curve. 

\begin{figure}[h]
\begin{center}
 \includegraphics[width=11cm]{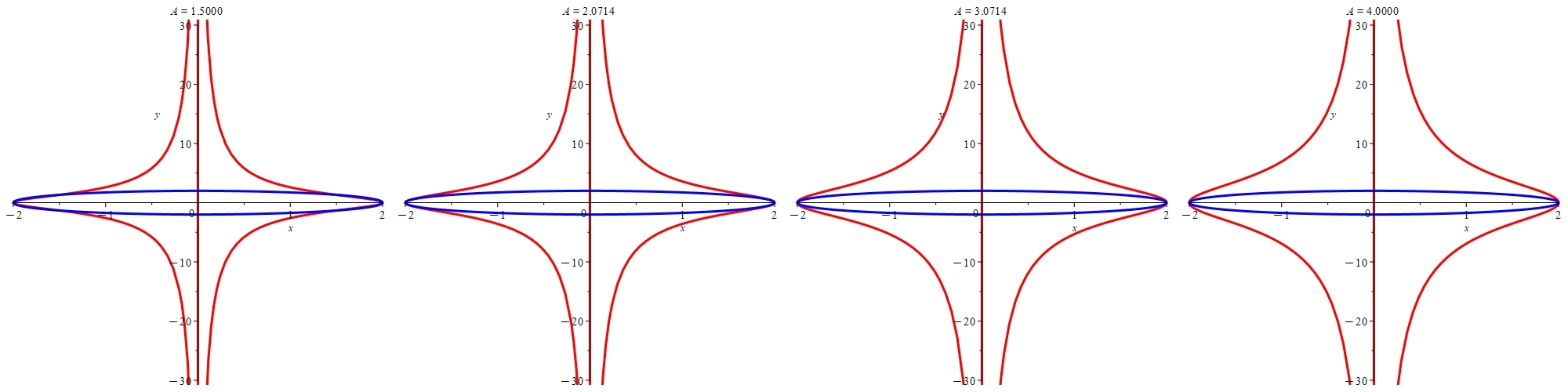}
 \caption{4 screenshots of the Maple animation}
\label{fig screenshots of the animation}
\end{center}
\end{figure}

Their differences are emphasized in Figure \ref{fig screenshots of the animation}, in an orthogonal but not orthonormal system of coordinates (therefore, the circle does not "look like" a circle) . These are screenshots of an animation programmed with Maple, with the \textbf{animate} command, as follows:
\small\begin{verbatim}
ac := animate(plot, [[r*cos(t), subs(d = A, yP), t = 0 .. 2*Pi]], 
       A = 0.5 .. 4, frames = 50, color = red, thickness = 3)
\end{verbatim} 

Here too, as in subsection \ref{subsection Hyperbolism of an ellipse}, the animation provides a visualization of the fact that for different values of the parameters, the obtained curve has the same topology.
 
\subsection{The circle passes through the origin}
\begin{figure}[htb]
\begin{center}
 \includegraphics[width=11cm]{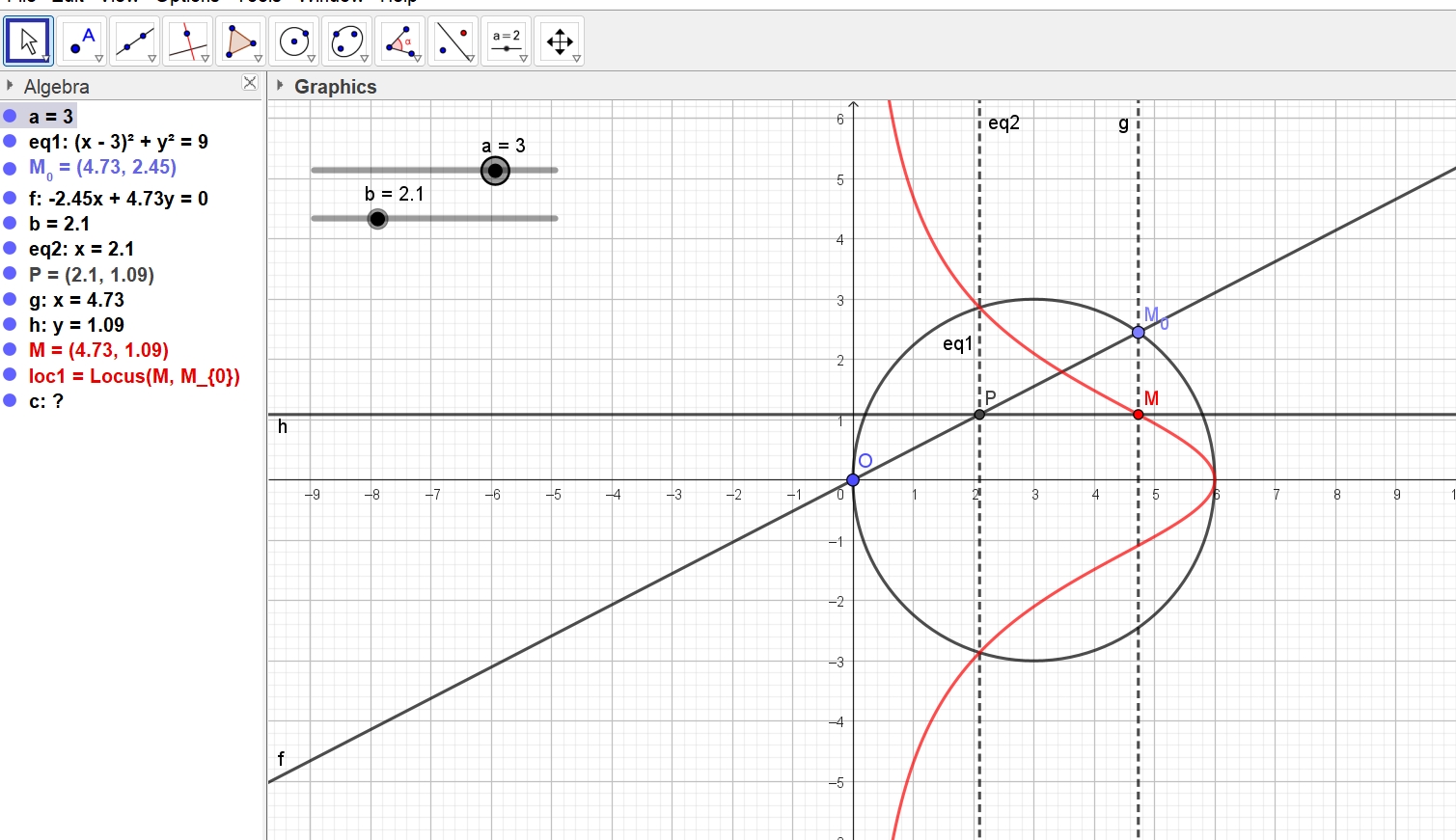}
 \caption{The line is a secant of the circle}
\label{fig line secant to the circle}
\end{center}
\end{figure}

We consider now a circle $\mathcal{C}$ passing through the origin and centered on the $x-$axis. The circle $\mathcal{C}$ has equation $\left(x-\frac {a}{2} \right)^2+y^2=\left(\frac{a}{2} \right)^2$. We take a line $D:x=b$, where $b<a$ and the point $O$ is the origin. The construction is displayed in Figure \ref{fig line secant to the circle}. The segments on the left are used instead of sliders.

A plot of the geometric locus of the point $M$ (the Tracer) when $M_0$ (the Mover) runs over the circle is obtained by the \textbf{Locus} command. Not as in  the previous section, the command \textbf{LocusEquation} does not provide an answer. Therefore, algebraic computations have to be performed using the CAS.

First, we derive a parametric presentation for the circle $\mathcal{C}$, using the intersection of the circle and of a line through the origin with slope $t$ (pay attention that this will be thw generic line $OM_0$)
\begin{equation}
\label{param eq for the circle}
\begin{cases}
x=\frac{a}{t^2+1}\\
y=\frac{at}{t^2+1}
\end{cases}
 \end{equation}

For the point $P$: if  $x_P=b$, then $y_P=tb$ and $M$ has coordinates 
\begin{equation}
\label{param M}
x=\frac{a}{t^2+1}\\
y=tb
\end{equation}
We define polynomials
\begin{equation}
\label{2polynomials}
P_1(x,y)=x(t^2 + 1) - a \qquad \text{and} \qquad P_2(x,y)=-b*t + y.
\end{equation}
Let $J=<P_1,P_2>$. By elimination of the parameter $t$, we obtain the following implicit equation of degree 3: 
\begin{equation}
\label{cubic equation}
xy^2+ b^2x-ab^2=0.
\end{equation}
Note that for $b=a$, i.e. the line is tangent to the circle at the end point of a diameter passing by the origin, the curve is a Witch of Agnesi.

\section{Piriform quartics}
\label{section piriform}

\subsection{Hyperbolism of a piriform quartic curve}
\begin{figure}[htb]
\begin{center}
 \includegraphics[width=11cm]{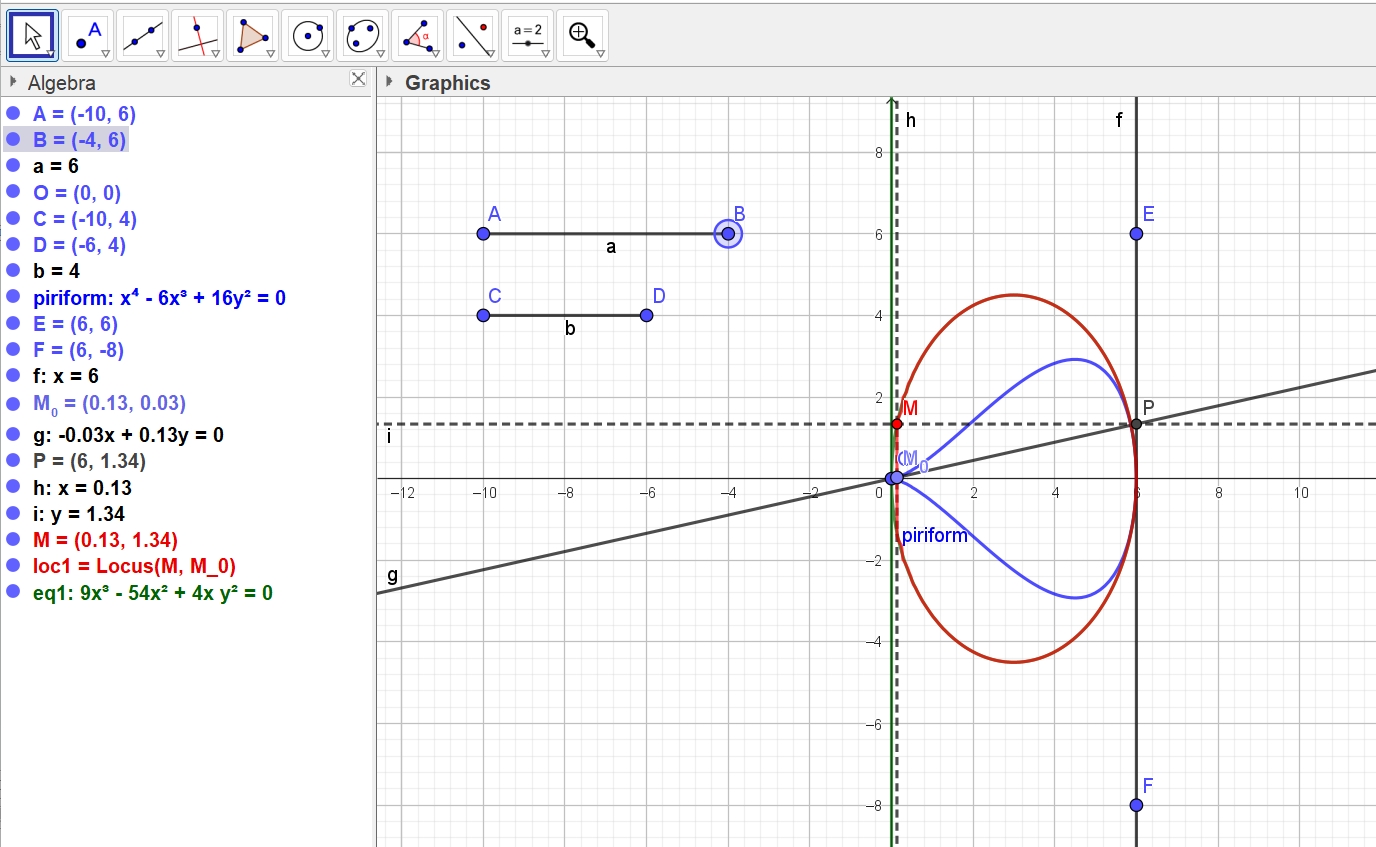}
 \caption{Hyperbolism of a piriform curve}
\label{fig hyperbolism of a piriform curve}
\end{center}
\end{figure}

We refer to \cite{mathcurve piriform} for the implicit and parametric presentation of the curve. A piriform curve is a quartic $\mathcal{P}$ whose equation is 
\begin{equation}
\label{eq piriform quartic}
b^2y^2=x^3(a-x),
\end{equation}
where $a$ and $b$ are positive parameters. In Figure \ref{fig hyperbolism of a piriform curve} (a screenshot of a GeoGebra applet\footnote{\url{https://www.geogebra.org/m/smcehrzy}}), these parameters are determined by 2 segments, in order to have a purely geometric construction.
This figure correspond to the case $(a,b)=(3,2)$. A plot is obtained and, simultaneously, an implicit equation, which reads as follows:
\begin{equation}
\label{eq hyperbolism of a piriform}
x(9x^2 - 54x + 4 y^2) = 0
\end{equation}
The factor $x$ determines the $y-$axis, emphasized in the Figure. This component is irrelevant to the geometric question; it appears because issues related to Zariski topology. Such issues are discussed in \cite{DPK-isoptics}, and are beyond the scope of the present article. The 2nd factor determines an ellipse, whose equation can be written as follows:
\begin{equation}
\label{eq ellipse_hyperbolism of a piriform}
\frac{(x-3)^2}{9}+\frac{y^2}{\frac{81}{4}}=1
\end{equation}
enabling to find the geometric characteristic elements of the ellipse. This is the desired hyperbolism of the piriform curve.

In order to perform the algebraic computations, it is worth to begin with a parametric presentation of the curve $\mathcal{P}$. We may use a trigonometric parametrization\footnote{It is given at \url{https://mathcurve.com/courbes2d.gb/piriforme/piriforme.shtml}.}:
\begin{equation}
\label{param eq piriform}
\begin{cases}
x= \frac{a}{2} (1+ \cos t) \\
y=\frac{a^2}{8b} (\sin 2t + 2 \sin t)
\end{cases}
\end{equation}
but we have to transform it into a rational parametrization. We choose a different way. Any line but the $y-$axis through the origin intersects again the curve $\mathcal{P}$. 
Let $L$ be the line whose equation is $y=tx, \; t \in \mathbb{R}$. We use the following Maple code:
\small
\begin{verbatim}
pear := b^2*y^2 - x^3*(a - x);
l := -t*x + y;
solve({l = 0, pear = 0}, {x, y});
par := allvalues(%[2]);
\end{verbatim}
\normalsize
The output gives two components, given by:
\small
\begin{equation}
\label{param intersection line and pear 1}
(x,y)=\left( \frac{a}{2}+\frac{\sqrt{-4 b^{2} t^{2}+a^{2}}}{2},\; t \left(\frac{a}{2}+\frac{\sqrt{-4 b^{2} t^{2}+a^{2}}}{2}\right) \right)
\end{equation}
\normalsize
and 
\small
\begin{equation}
(x,y)=\left( \frac{a}{2}-\frac{\sqrt{-4 b^{2} t^{2}+a^{2}}}{2}, t \left(\frac{a}{2}-\frac{\sqrt{-4 b^{2} t^{2}+a^{2}}}{2}\right)\right)
\end{equation}
\normalsize
We work now with the first component, the 2nd one can be treated exactly in the same way. 

Denote $x_P=a, \; y_P=ta$. By construction, we have:
\begin{equation}
\label{coordinates M}
\begin{cases}
x_M=\frac{a}{2} + \frac {1}{2} \sqrt{-4b^2t^2 + a^2}\\
y_M=ta
\end{cases}
\end{equation}
Now we use the following Maple code:
\small
\begin{verbatim}
p1 := x - a/2 =sqrt(-4*b^2*t^2 + a^2)/2;
p1 := p1^2;
p1 := lhs(p1) - rhs(p1);
p2 := -a*t + y;
J:=<p1,p2>;
JE := EliminationIdeal(J, {a, b, x, y})
G := Generators(JE)[1];
\end{verbatim}
\normalsize
whose output is the polynomial
\begin{equation}
G(x,y) =  a^2x^2 -a^3x + b^2y^2
\end{equation}
We have 
\begin{equation}
G(x,y)=a^2(x^2-ax)+b^2y^2=a^2 \left( x-\frac {a}{2} \right)^2 +b^2y^2-\frac{a^4}{4}.
\end{equation}
Rewriting the equation under the form
\begin{equation}
\frac {\left( x-\frac {a}{2} \right)^2}{b^2}+\frac{y^2}{a^2}=\frac{a^2}{4b^2},
\end{equation}
which is equivalent to
\begin{equation}
\frac{\left( x-\frac{a}{2} \right)^2}{\left(\frac{a}{2}\right)^2}+\frac{y^2}{(2b)^2}=1,
\end{equation}
we identify that the requested curve is an ellipse.

\begin{remark}
The obtained hyperbolism has been identified as an ellipse by algebraic means. This appears also in the GeoGebra applet, but a slight modification may induce a big change, and the implicit equation may not be obtained. An experimental way to check that the curve is an ellipse consists in marking 5 arbitrary points on the curve (with the \textbf{Point on Object} command) and  determine a conic by 5 points on it (there is a button driven command for this). This is a numerical checking, not a symbolic proof, but may enable students to proceed further. Of course, the algebraic computations that we performed with the CAS are a must.     
\end{remark}

\subsection{A piriform curve as antihyperbolism of a circle}
\label{subsection a different construction leading to a piriform}
The choice of the original curve in previous subsection incites to have a look at piriform curves from another point of view; see \cite{hyperbolism mathcurve,mathcurve piriform}.
 \begin{definition}
 \label{def antihyperbolism}
With the notations of Definition \ref{def hyperbolism}, the inverse transformation $ (x,y) \rightarrow \left( x, \frac{xy}{a} \right)$ is called \emph{antihyperbolism}.
 \end{definition}

\begin{figure}[htb]
\begin{center}
 \includegraphics[width=12cm]{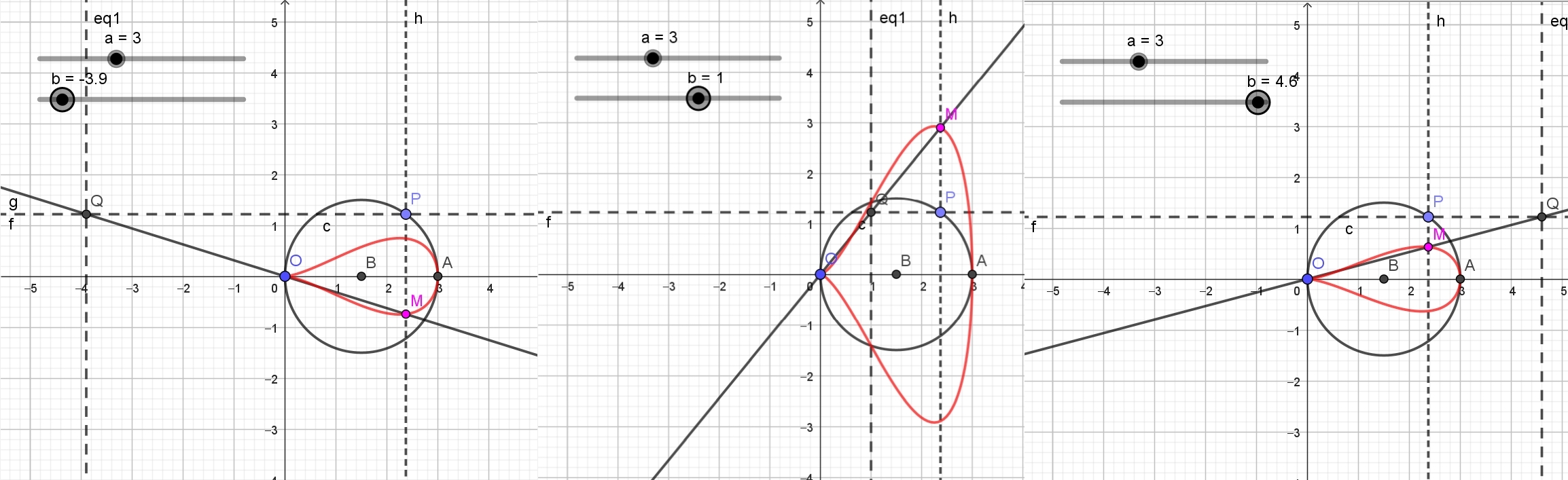}
 \caption{A piriform quartic constructed from a circle - exploration with the sliders}
\label{fig pear 3 screenshots}
\end{center}
\end{figure}

Screenshots of a GeoGebra applet\footnote{\url{https://www.geogebra.org/m/yfq7wmqz}} are displayed in Figure \ref{fig pear 3 screenshots}.
Exploration using the sliders show that for $|b|> a$, the piriform curve lies inside the circle and is tangent to it at $A$. For $-a<b<a$, the curve is part inside and part outside of the circle, and is still tangent to the circle at $A$, but from outside.  

\section{Hyperbolism of a nephroid}

\begin{figure}[htb]
\begin{center}
\includegraphics[width=10cm]{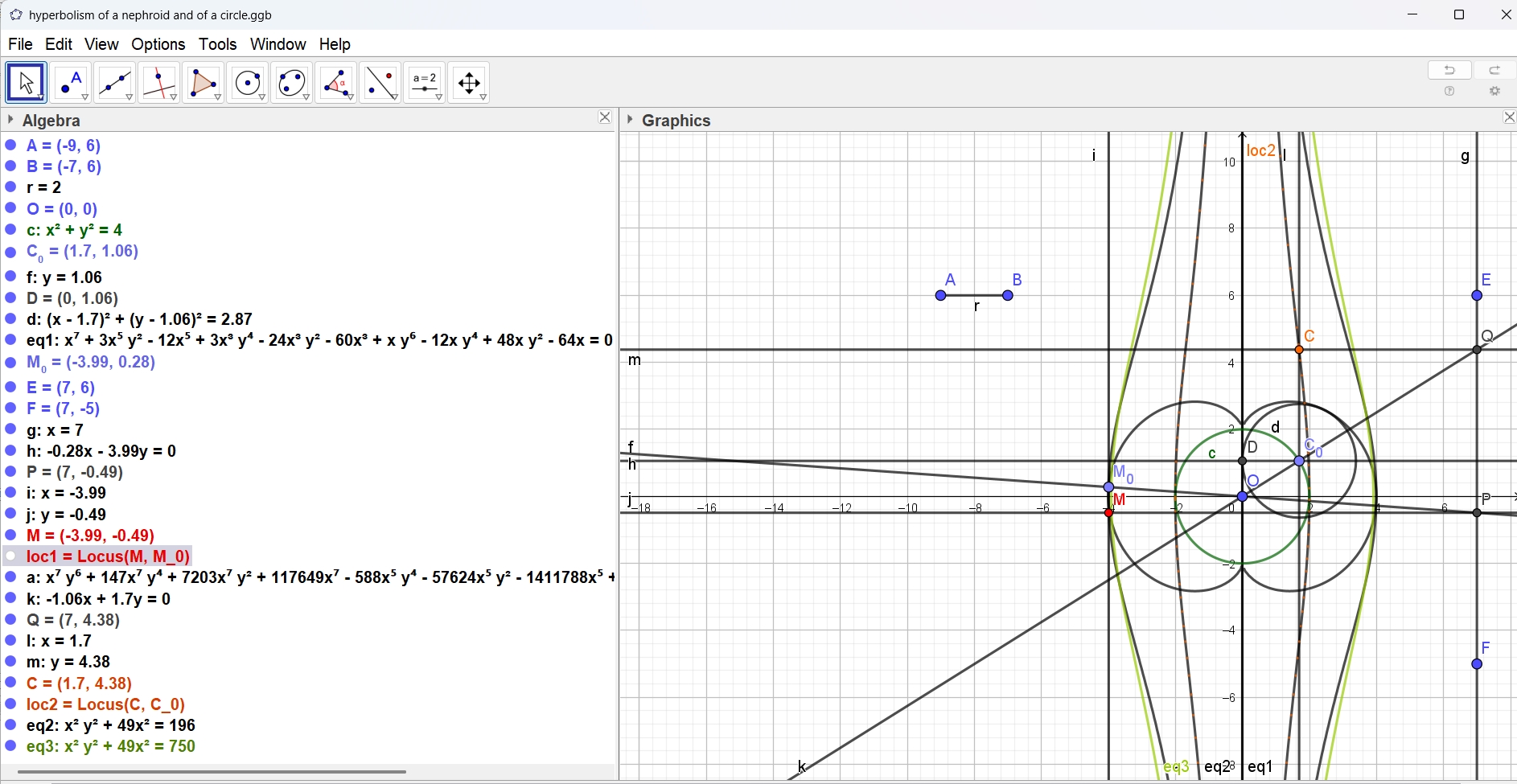}
\caption{Comparison  between hyperbolisms of a nephroid and of a circle}
\label{fig hyperbolism of a nephroid}
\end{center}
\end{figure}

Let $\mathcal{C}$ be a circle whose center is at the origin. The envelope of the family of circles centered on $\mathcal{C}$ and tangent to the $y-$ axis is a sextic called a nephroid \cite{revival}. By that way the curve has been constructed in Figure \ref{fig hyperbolism of a nephroid}, which is a screenshot of a GeoGebra applet\footnote{\url{https://www.geogebra.org/m/pyhk9qvr}}. 
We denote the nephroid by $\mathcal{N}$. 
Later, this figure will enable to compare the hyperbolism of  $\mathcal{N}$ with a hyperbolism of a circle, as described in Subsection \ref{subsection circle 1}.

\begin{figure}[htb]
\begin{center}
\includegraphics[height=5cm]{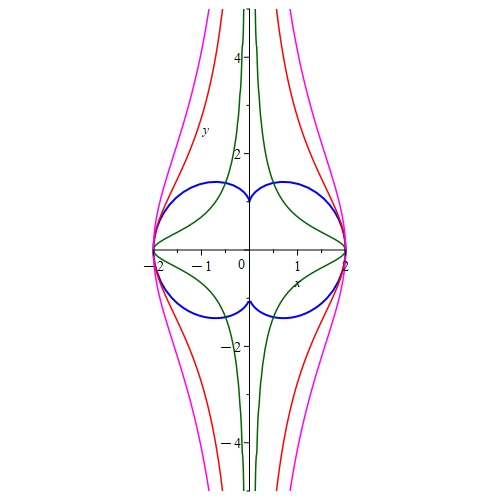}
\caption{Three hyperbolisms of the same nephroid with respect to the origin}
\label{fig hyper nephroid}
\end{center}
\end{figure}

A general implicit equation for a nephroid is
\begin{equation}
\label{implicit eq nephroid}
4(x^2+y^2-a^2)^3=27a^4x^2,
\end{equation}
where $a$ is a positive parameter. The curve $\mathcal{N}$ can also be described by a parametric presentation
\begin{equation}
\label{param eq nephroid}
\begin{cases}
x= 2a \sin^3 t \\
y= a(1+\sin^2 t)\; \cos t
\end{cases}
\end{equation}
To construct a hyperbolism of $\mathcal{N}$ with respect to the origin and to a vertical line, we follow the same path as in previous sections. This begins with deriving from Equation (\ref{param eq nephroid}) a rational parametrization. Applying the substitution in Equation (\ref{rational substitution for trigo}), we obtain:
\begin{equation}
\label{rational param nephroid}
\begin{cases}
x_0=\frac{-u^6 - 9u^4 + 9u^2 + 1}{(u^2 + 1)^3}\\
y_0=\frac{16au^3}{(u^2 + 1)^3}
\end{cases}
\end{equation}
The coordinates of $P$ are thus
\begin{equation*}
\begin{cases}
x_P= b\\
y_P= \frac{(-u^6 - 9u^4 + 9u^2 + 1)b}{16au^3}
\end{cases}
\end{equation*}
and the coordinates of $M$ are
\begin{equation}
\label{param eq hyperbolism nephroid}
\begin{cases}
x_M= \frac{16au^3}{u^2 + 1)^3}\\
y_M= \frac{(-u^6 - 9u^4 + 9u^2 + 1)b}{16au^3}
\end{cases}
\end{equation}
We run now the following Maple code, similar to what we did already:
\small
\begin{verbatim}
p1 := x*denom(xM) - numer(xM);
p2 := y*denom(yM) - numer(yM);
J := <p1, p2>;
JE := EliminationIdeal(J, {a, b, x, y});
G := Generators(JE)[1];
evala(AFactor(G));
\end{verbatim}
The last command is intended to check that the obtained polynomial $G(x,y)$ is irreducible. We have here a polynomial of degree 12, which fits Remark \ref{remark degree}:
\begin{equation}
\label{implicit eq hyperbolism  of nephroid}
\begin{split}
G(x,y)=4a^6x^6y^6 - 12a^6b^2x^4y^4 + 12a^4b^2x^6y^4 + 12a^6b^4x^2y^2 - 24a^4b^4x^4y^2 \\
+ 12a^2b^4x^6y^2 - 4a^6b^6 - 15a^4b^6x^2 - 12a^2b^6x^4 + 4b^6x^6. \qquad \qquad
\end{split}
\end{equation}
Figure \ref{fig hyper nephroid} shows the curve $\mathcal{N}$ for $a=1$ and hyperbolisms with respect to the origin and the line whose equation is $x=b$ for $b=1/2, 2,3$.

\section{A construction of a lemniscate}
\label{lemniscate}
The construction proposed in this section is slightly different from a hyperbolism. Studying at least one case may be an appeal to explore more cases.
   
\begin{figure}[htb]
\begin{center}
 \includegraphics[width=10cm]{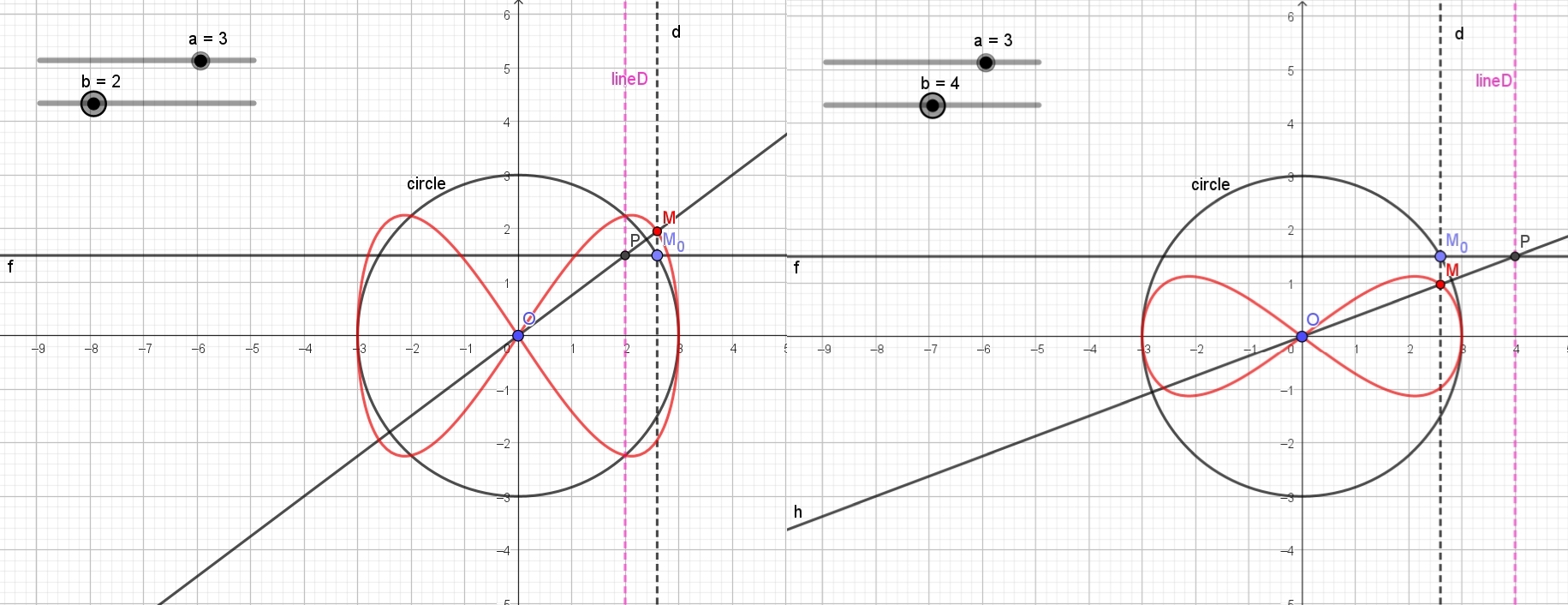}
 \caption{The circle is centered at the origin and the line intersects it.}
\label{fig center at the origin - line secant to the circle}
\end{center}
\end{figure}

We consider a circle $\mathcal{C}$  whose center is the origin and radius $r$ and a line $D$ whose equation is $x=b$. Take a point $M_0$ on the circle $\mathcal{C}$; the point $P$ is the intersection of $D$ with the horizontal line through $M_{0}$. Then we define $M$ to be the point of intersection of the line $(OP)$ with the vertical line through $M_0$. The geometric locus of $M$ when $M_0$ runs on  $\mathcal{C}$ is a lemniscate. This is illustrated in Figure \ref{fig center at the origin - line secant to the circle}. The figure shows screenshots of a GeoGebra applet\footnote{\url{https://www.geogebra.org/m/zkgmxume}}; using the sliders may help to study the family of curves.

For $b<r$, the lemniscate lies part out of the circle. For $b>r$, the lemniscate is bounded by the circle. In both cases, the lemniscate is tangent to the circle at its vertices, i.e. its points on the $x-$axis different form the origin.

GeoGebra could not give an implicit equation for the lemniscate, we will look for such an equation via algebraic computations. We choose the same parametrization as above for the circle, i.e. the coordinates of $M_0$ are given by:  
\begin{equation}
\label{param circle}
(x_0,y_0)=\left( a \cdot \frac{1-u^2}{1+u^2},\; a\cdot \frac {2u}{1+u^2} \right).
\end{equation}
Then the coordinates of $P$ are
\begin{equation*}
(x_P,y_P)=\left( b, a \cdot \frac {2u}{1+u^2} \right)
\end{equation*}
The equation of the line $(OP)$ is $y=\frac{2au}{b(1+u^2)} \cdot x$, whence the coordinates of $M$: 
\begin{equation}
\label{coord M - lemniscate}
(x_M,y_M)=\left(a \cdot \frac{1-u^2}{1+u^2}, \; \frac {2a^2u(1-u^2)}{b(1+u^2)^2} \right).
\end{equation}
We apply the following code:
\small
\begin{verbatim}
p1 := x*denom(xM) - numer(xM);
p2 := y*denom(yM) - numer(yM);
J := <p1, p2>;
JE := EliminationIdeal(J, {a, b, x, y});
G := Generators(JE)[1];
\end{verbatim}
\normalsize
and obtain the polynomial
\begin{equation}
\label{eq lemniscate}
G(x,y)=a^2x^4-a^4x^2 + a^2b^2y^2. 
\end{equation}
The vanishing set of $G(x,y)$ in the real plane is called a \emph{lemniscate of Gerono}.

\section{Special features of the curves which have been obtained}

In Section \ref{section circle and ellipse}, the obtained hyperbolisms show common features: 
\begin{enumerate}
\item The curve has two disjoint components. On the one hand, these components cannot be distinguished by algebraic means as the defining polynomial is irreducible over the field of real numbers. This can be proven using Maple command \textbf{evala(AFactor(...)}; the need for that command, and not the ordinary \textbf{Factor} command, has been analyzed in \cite{DPR}. 
\item The geometric basis of the construction induces easily that the two components are symmetric about the $y-$axis.
\item From the exploration, the $y-$axis seems to be an asymptote to both. Actually, this is a consequence of the construction: when the point $M_0$ gets arbitrarily close to the $y-$axis, the slope of the line $(OM_0)$line, and thus  $y-$coordinate of $M$, tend to infinity.
\end{enumerate}

Actually the origin does not belong to the hyperbolism of the original curve. It can be obtained only if $M_0$ is at the origin, but in this case the line $(OM_0)$ is not defined.

 \section{Discussion}

The 2 first examples may be treated by hand, with students having learnt an elementary course in the spirit of \cite{adams and loustaunau}.

In this paper, we considered a question in plane geometry, which can look rather simple. Its translation using technology revealed a more complicated situation that foreseen.  On the other hand, the mathematics "behind" the screen are heavier than the general domain announced. 

The problem can be understood by a regular High-School student\footnote{Decades ago, such constructions were performed in High-School by hand using paper, pencil and ruler. Only in easy cases, an equation was derived.}, but the algebraic machinery beneath is then out of reach, as it belongs to Computer Algebra, and uses algorithms from Ring Theory.  Moreover, we deal here with curves of higher degree than what High-School students learn, as explained in section \ref{intro}.

Along the paper, implicitization has been performed using polynomial rings and elimination. Maple has an \textbf{implicitize} command, but we preferred to have our computations more understandable, and not to use this command as a blackbox. The user has no access to the algorithms, but via reading the paper in reference on the help page of the command. One the one hand, an educator has to make a decision between low-level and high level commands (see \cite{low level}). On the other hand, what can be done depends on the students' theoretical background; sometimes it is necessary to use the CAS in order to bypass a lack of knowledge (see \cite{bypass}). 

Moreover, respective affordances of the two kinds of software that we used have been analyzed in the past, and the importance of networking between them has been emphasized \cite{eugenio,eugenio 2,DPK-networking}. See also \cite{weinhandl et al} for the more recent developments with GeoGebra.  Here the strengths of both kinds appeared sometimes different from what we were accustomed. Anyway, it is only natural that for different questions, the respective abilities and strengths of the different kinds of software which are used may vary. 

Finally, the lemniscate of Gerono in Section \ref{lemniscate} shows that classical constructions can be imitated with slight modifications and lead to software's broader experience. After all, curiosity is the main engine for exploration and discovery.    

\section*{References}

\end{document}